\documentclass{amsart}
\usepackage{amssymb}

\usepackage{hyperref}

\newcommand{\Z}{\mathbb Z}
\newcommand{\Q}{\mathbb Q}

\def\as#1{\renewcommand\arraystretch{#1}}

\def\diso{\lower.4ex\hbox{$\downarrow$}\raise.4ex\hbox{\mbox{\scriptsize
$\wr$}}}

\def\ff#1{\mathbb{F}_{#1}}

\def\iso{\ \lower.3ex\hbox{\as{.08}$\begin{array}{c}\lra\\\mbox{\tiny $\sim\,$}\end{array}$}\ }

\def\lg{l\raise.6ex\hbox to.2em{\hss.\hss}l}
\def\lra{\longrightarrow}

\def\orb{\hbox to  .3em{$\backslash$}\backslash}

\def\p{\mathfrak{p}}

\def\ty{\mathbf{t}}

\newcounter{cs}
\stepcounter{cs}
\newcommand{\casos}{\begin{itemize}}
\newcommand{\fcasos}{\end{itemize}\setcounter{cs}{1}}

\newfont{\tit}{cmr12 scaled \magstep3}

\def\pot{\textasciicircum}

\def\magma#1{\vskip 3mm\hskip.6cm{\tt \begin{minipage}{11.6cm}  #1  \end{minipage}}\vskip 3mm}

\def\magmad#1{\vskip 3mm\hskip3mm{\tt \begin{minipage}{15cm}  #1  \end{minipage}}\vskip 3mm}

\begin{document}
\title{Arithmetic in big number fields:\\ the {\tt '+Ideals'} package}
\author[Gu\`ardia]{Jordi Gu\`ardia}
\address{Departament de Matem\`atica Aplicada IV, Escola Polit\`ecnica Superior
d'Enginyera de Vilanova i la
Geltr\'u, Av. V\'\i ctor Balaguer s/n. E-08800 Vilanova i la Geltr\'u,
Catalonia}
\email{guardia@ma4.upc.edu}

\author[Montes]{\hbox{Jes\'us Montes}}
\address{Departament de Ci\`encies Econ\`omiques i Socials,
Facultat de Ci\`encies Socials,
Universitat Abat Oliba CEU,
Bellesguard 30, E-08022 Barcelona, Catalonia, Spain}
\email{montes3@uao.es}

\author[Nart]{\hbox{Enric Nart}}
\address{Departament de Matem\`{a}tiques,
         Universitat Aut\`{o}noma de Barcelona,
         Edifici C, E-08193 Bellaterra, Barcelona, Catalonia, Spain}
\email{nart@mat.uab.cat}
\thanks{Partially supported by MTM2009-13060-C02-02 and MTM2009-10359 from the
Spanish MEC}
\date{}

\makeatletter
\@namedef{subjclassname@2010}{%
  \textup{2010} Mathematics Subject Classification}


\begin{abstract}
We introduce our package \texttt{+Ideals} for Magma, designed to perform the basic tasks related to  ideals in number fields without pre-computing integral bases. It is based on Montes algorithm and a number of local techniques that we have developed  in a series of papers in the last years.
\end{abstract}

\maketitle

%\tableofcontents

\section*{Introduction}

A commonplace in Number Theory problems is the need of effective computations. The last decade has seen an exponential growth of mathematical software to satisfy this necessity. Striking ideas have expanded our computational limits far beyond our present requirements in many fields, specially those related to cryptographic applications. However, there remain some parts of Computational Number Theory which have not evolved that much. Algebraic Number Theory appears to be one of these.

While the available algorithms for working over number fields are bright and efficient, the computational requirements of {\rm real life} problems in Number Theory frequently exceed their capability. As soon as we have to work in a number field defined by a polynomial with large discriminant, the determination of its  ring of integers, which is the cornerstone of all the {\rm classic} algorithms, becomes unfeasible, thus making impossible further computations.

We have developed a number of techniques, based on higher order Newton polygons, to deal with the prime ideals in number fields  (\cite{GMNalgorithm}, \cite{GMNbasis}, \cite{GMNokutsu}, \cite{GMN3}).
The application of these ideas
leads to extremely fast and efficient algorithms for the basic tasks in Computational Number Theory.  These new algorithms can work  over number fields of degree 1000 or with a 300-digit discriminant with a home computer in a few seconds. A new philosophy is beyond this project: the centre of attention in Computational Algebraic Number Theory may be not the ring of integers, but the prime ideals.

We have implemented these algorithms in the {\tt +Ideals} package for Magma. Besides factorization of ideals, %and determination of the discriminant, it
this package gives facilities to compute valuations at prime ideals, reduction modulo prime ideals and Chinese Remainder Theorem problems. The algorithms have been presented in \cite{GMN3}. The package can also compute integral bases of number fields, %under a certain conjecture,
following the algorithm introduced in \cite{GMNbasis}. The goal of this paper is to introduce the package, presenting some examples of use and describing its main functions.

The paper is structured in three sections: after giving
a general overview of the package in the first section, we illustrate its use and capabilities by means of some concrete examples in Section 2. Section 3 gives a more exhaustive description of routines of the package.

A primitive version of the package was published about two years ago, under the name {\tt Newton}. It included essentially the decomposition of primes in number fields and the computation of integral bases. The functions in the \texttt{Newton} package remain in the {\tt +Ideals} package, to ensure compatibility for all users.

The \texttt{+Ideals} package can be obtained on demand from the authors, or downloaded from its web page
 ({\url  {http://www-ma4.upc.edu/~guardia/+Ideals.html}}). The package requires version v2.15 of Magma or newer.

\section{General overview of the package}

The {\tt {+Ideals}} package is concerned with the arithmetic of fractional ideals in number fields. The keystone  is Montes algorithm, which provides local information for the decomposition of rational primes in number fields. The great advantage of the algorithm is to avoid the use of integral bases. This shortcut makes possible to perform fast computations in number fields of huge degree and index.

The underlying philosophy of the package is that prime ideals are the basic tools for many computations in  number fields, and hence all the functions in the package perform tasks related to them. The main task, of course, is factorization of fractional ideals, but computation of valuations and reduction maps are also important. As soon as a prime ideal is detected along a computation, the information which can be useful for other computations is stored in some new attributes for number fields, so that it can be easily retrieved if needed.

The package introduces two new data types, {\tt PrimeIdealRecord} and {\tt Ideal\-Record}, designed to represent  fractional ideals in number fields, and some attributes for number fields to store information about these ideals. Most of the functions included in the package are extended versions of the usual Magma routines for ideal arithmetic,
designed to deal with these new data types, so that the experienced user of Magma should find no difficulty in the use of the package.

Let us finish this overview with a self-explanatory small piece of code illustrating how to perform the most common tasks in ideal arithmetic with our package:
\magma{
> Zx<x>:=PolynomialRing(Integers());\\
> K:=NumberField(x\pot100-x\pot75+x\pot50+2\pot500);\\
> I:=ideal(K,2);\\
> time Factorization($\sim$I);       \\
Time: 1.300 \\
>I\`{}FactorizationString;     \\
P(2,1)\pot2*P(2,2)\pot2*P(2,3)\pot2*P(2,4)*P(2,5)*P(2,6)*P(2,7)*P(2,8)\\
> K\`{}LocalIndex[2];\\
12250
}

\newpage
\section{Some hard examples}

We illustrate the use of the package by means of some examples picked from different topics in Computational Number Theory. The computations have been made with Magma v2.15-11  in a Linux server, with two Intel Quad Core processors,
running at 3.0 Ghz, with 32Gb of RAM memory. All the running times mentioned in these examples are expressed in seconds.

\subsection{A {\em real life} problem}

A very common problem in Arithmetic Geometry is the determination of good models for curves over number fields. Usually, a curve with potentially good reduction acquires good reduction over a number field generated by certain torsion points in the jacobian of the curve. The number fields encountered in this process are usually huge, and it is a difficult computational problem to find the places over the rational primes we are interested in. We show in this example how to use  our package to solve such a problem. We will determine the reduction type of a rational elliptic curve over a huge number field, applying Tate's algorithm. The decomposition of primes in the number field is the minimal information required, but no integral basis is necessary to run the algorithm, so that the problem fits very well  the philosophy of our package.

Consider the  elliptic curve
$$
E: y^2=x^3+3x^2+3x,
$$
and let $K=\Q(\theta)$ be the number field determined by a root $\theta\in\overline{\Q}$ of the 17-th division polynomial. It is a number field of degree 144. We can determine it with the following Magma code:

\magma{
> E:=EllipticCurve([0,3,0,3,0]);\\
> f:=DivisionPolynomial(E,17);\\
> K<w>:=NumberField(f);
}

The discriminant of the polynomial $f$ defining $K$ is huge, but easy to factor, due to its geometric origin:

\magma{
> d:=Discriminant(f); \newline
> Factorization(d); \newline
 [ <2, 13728>, <3, 10296>, <17, 143> ]
}

We cannot expect to find the ring of integers $\Z_K$ of such a big number field with the standard Magma routines: 
 the function {\tt MaximalOrder } exhausts the computer's memory after some hours of work. But we do can compute in this ring! We must redefine $K$ through a monic polynomial with integer coefficients; since  the leading coefficient of $f$ is 17, we simply type:

\magma{
> Qx<x>:=PolynomialRing(Rationals());\\
> g:=17\pot143*Evaluate(f,x/17);\\
> K<w>:=NumberField(PolynomialRing(Integers())!g);
}
\noindent
Note that  a root of $g$ is $w=17\theta$.
Let us now attach our package\footnote{Replace {\tt dir} with the name of the directory where you saved the package.}:

\magma{ > Attach("dir/+ideals.m");}

\noindent
Let us find  the decomposition of 2 in $\Z_K$. We first define the ideal $2\Z_K$ using the new function {\tt ideal}, and then we factor it into prime ideals:

\magma{
> i:=ideal(K,2);\\
> time Factorization($\sim$i);\\
Time: 0.010\\
> i\`{}FactorizationString;\\
P(2,1)\pot3*P(2,2)\pot3*P(2,3)\pot3*P(2,4)\pot3*P(2,5)\pot3*P(2,6)\pot3\\
}
\noindent
The result means that $2\Z_K$ is the cube of the product of six different prime ideals: $2\Z_K=(\p_1\cdots\p_6)^3$.
The use of $\sim$ allows the factorization routine to append information to the record corresponding to the given ideal. A shorter way to find the decomposition of a prime in $\Z_K$ is:

\magma{
 > time Factorization(ideal(K,3)); \newline
[
    [ 3, 1, 2 ],
    [ 3, 2, 2 ],
    [ 3, 3, 2 ],
    [ 3, 4, 2 ],
    [ 3, 5, 2 ],
] \newline
Time: 0.010
}
\noindent
which indicates that there are  5 different prime ideals in $\Z_K$ over 3, all of them with ramification index 2 (but now we can not access the assigned attributes of the ideal, since they are not saved)
We can access these prime ideals through the new attribute {\tt PrimeIdeals} established for number fields. The command:

\magma{
> K\`{}PrimeIdeals[2];
}
\noindent
will generate a long output with all the data already computed for these ideals. For instance:

\magma{
> K\`{}PrimeIdeals[2,1]`f;\\
8
}
\noindent
gives the residual degree of the first of the prime ideals dividing 2. The residual degrees of the prime ideals above 3 are:
\magma{
> [P\`{}f:P in K`PrimeIdeals[3]];\newline
[ 8, 16, 16, 16, 16 ] \newline
}
\noindent
We already have generators for the ideals at hand; typing:
\magma{
> alpha:=K\`{}PrimeIdeals[3,1]\`{}Generator;alpha;\\
1/3433683820292512484657849089281*(w\pot129+2451*w\pot128+...
}
\noindent
we  obtain an element $\alpha\in\Z_K$ such that $\p=(3,\alpha)$ is the prime over 3 with residual degree 8.

Along the factorization routine of the ideal $p\Z_K$, the program attaches to each prime ideal over $p$ a local \emph{type}, which carries on some canonical invariants  (\cite{GMNokutsu}). These types are stored as an attribute of the prime ideals ({\tt K\`{}PrimeIdeals[p,k]\`{}Type}) and contain essential information.
As a by-product of the computation of these types we obtain the $p$-index of the polynomial $g$. It is saved in an associative array of Magma called {\tt LocalIndex}, with one entry for each prime number already factored in $\Z_K$:
\magma{
>K\`{}LocalIndex[2];\\
6816
}
\noindent
Using this array, we can easily find the  discriminant of $K$; we first factor completely the discriminant $d$ of $g$ and then get the global index:
\magma{
> time Factorization(ideal(K,17));\\
...\\
Time: 0.380\\
> ind:=(2\pot (K\`{}LocalIndex[2])*3\pot(K\`{}LocalIndex[3])*\\
\qquad\qquad17\pot(K\`{}LocalIndex[17]));\\
> disc:=d div ind\pot2; \\
> Factorization(disc);\\\,
>  [ <2, 96>, <3, 72>, <17, 143>] \\
}

We can check that the point $P=(w/17,\sqrt{w^3/17^3+3w^2/17^2+3w/17})\in E(\overline{K})$ is not defined over $K$. The computation of the square root in the $y$-coordinate may be quite expensive, so that we simply reduce its second coordinate modulo a prime ideal:
\magma{
> Factorization(ideal(K,5));\\
> yred:=(w\pot3/17\pot3+3*w\pot2/17\pot2+3*w/17) mod K\`{}PrimeIdeals[5,1];\\
> IsSquare(yred);\\
false
}

Let us apply now Tate's algorithm  to determine the reduction type of $E$ over the  prime $\p=(3,\alpha)$ defined above. We have $a_1=a_3=a_6=0$, $a_2=a_4=3$, $b_2=12$, $b_4=6$, $b_6=0$, $b_8=-9$. It's a simple matter to check the cumulative hypothesis of Tate's algorithm:
\magma{
> P:=K\`{}PrimeIdeals[3,1];\\
> dE:=Discriminant(E);\\
> PValuation(dE,P);\\
6 \newline
> [PValuation(ai,P): ai in aInvariants(E)];\newline
[ Infinity, 2, Infinity, 2, Infinity ]\newline
> [PValuation(bi,P): bi in bInvariants(E)];\newline
[ 2, 2, Infinity, 4 ]
}
\noindent
Thus, we have
$$
\begin{array}{l}
v_\p(\Delta)>0\\
v_\p(a_3)>0, v_\p(a_4)>0, v_\p(a_6)>0,\\
v_\p(b_2)>0,\\
v_\p(a_6)\ge 2,\\
v_\p(b_8)\ge3,\\
v_\p(b_6)\ge3,\\
v_\p(a_1)>0, v_\p(a_2)>0, v_\p(a_3)\ge2, v_\p(a_4)\ge2, v_\p(a_6)\ge3.\\
\end{array}
$$
We must now consider the polynomial
$$
P(T)=T^3+a_{2,1}T^2+a_{4,2}T+a_{6,3},
$$
where $a_{i,m}=a_i/\Pi^m$, and $\Pi$ is a local uniformizer for $\p$. We get one by typing:
\magma{
> PI:=P\`{}Generator;
}
\noindent
Since $a_{2,1}\equiv a_{6,3}\equiv0(\operatorname{mod} \p)$ and $v_\p(a_{4,2})=0$, the polynomial $P(T)(\,\operatorname{mod} \p)$ has three different roots in an algebraic closure of the residue field $k(\p)$ of $\p$, so we already know that the Kodaira symbol of $E$ is $I_0^\ast$. The local conductor is $\p^{v_\p(\Delta)-4}=\p^2$. To determine the local index $c_\p=[E(K_\p)/E_0(K_\p)]$, we need the number of roots of $P(T)$ in $k(\p)$. We may proceed as follows\footnote{The computational representation of the residue field of a prime ideal is computed as a successive extension of finite fields: $$\ff0\subseteq \ff1\subseteq \cdots\subseteq \ff{r+1}=\ff{\p}=\ff0[z_0,z_1,\dots,z_r],$$ where $\ff{0}$ is the prime field and $\ff{i}:=\ff{i-1}(z_{i-1})$, for $1\le i\le r+1$. The parameter $r$ is the \emph{order} of the type attached to $\p$.
}
:
\magma{
> time  m:=3/PI\pot2 mod P;      \\
Time: 13.570  \\
>  F:=ResidueField(P); \\
> FT<T>:=PolynomialRing(F); \\
> Factorization(T\pot3+m*T);\newline
[    <T, 1>, \newline
 <T\pot2 + 2*z1\pot7 + 2*z1\pot6 + 2*z1\pot4 + 2*z1\pot3 + z1 + 2, 1>
 ]
}
\noindent
Hence, the local index is $c_\p=2$. It is easy to check that for the remaining prime ideals  in $\Z_K$ above 3, the algorithm proceeds in the same way (all the ramification indices are 2) up to this point, where we find that the corresponding polynomial $P(T)$ has three different roots in the residue field, and hence the local indices at these primes are equal to 4.

Finally, our package has its own routine {\tt IntegralBasis}  to compute a $\Z$-basis of $\Z_K$:

\magma{
> time  basis:=IntegralBasis(K);      \\
Time: 18.620
}

It has to be mentioned, however, that this routine is based on a conjecture that we have not been able to prove yet (\cite{GMNbasis}). In any case, the routine always returns a correct result (it performs a check of its output) and is prepared to warn the user in the improbable case that the algorithm fails (after millions of tests, we have not encountered any counterexamples).

\subsection{A {\em tall} number field}\label{example2}

In this second example we will force the computational capabilities of the package to show its power. The concrete computations are not significant on its own; they have been chosen to illustrate the different routines of the package.

Consider the number field $K=\Q(\theta)$ given by a root  $\theta\in\overline\Q$ of the polynomial $f(x)=x^{1000} + 2^{50} x^{50} + 2^{60}$. The factorization of the discriminant of $f$ is
$$
\mbox{Disc}(f)=
2^{53940}3^{50}5^{2000}127^{50} 313^{50} 743^{50} 4886229527^{50}p^{50},
$$
with $p=337572698551220494882323528404563236947916489629537$.
The large degree of $f$ makes impossible to work in this number field using the standard functions of Magma, even after factorizing the discriminant, since the computation of the integral basis is compulsory for these functions. But our functions avoid this computation, so that we can work with ideals in $K$. Let us begin by determining the discriminant of $K$:
\magma{
> Attach("/dir/+Ideals.m");\\
> Z:=Integers();\\
> Zx<x>:=PolynomialRing(Integers());\\
> f := x\pot1000 + 2\pot50* x\pot50 + 2\pot60;\\
> d:=Discriminant(f);\\
> Factorization(d);\newline
[ <2, 53940>, <3, 50>, <5, 2000>, <127, 50>, <313, 50>, <743, 50>, <4886229527,
50>, <337572698551220494882323528404563236947916489629537, 50> ] \newline
> badprimes:=PrimeDivisors(d);\\
> K<w>:=NumberField(f);\\
> for p in badprimes do time Factorization(ideal(K,p)); end for;\\
...
}
\magma{
>[K\`{}LocalIndex[p]: p in K\`{}FactorizedPrimes];\newline
[ 26235, 0, 20, 0, 0, 0, 0, 0 ]
}
The last command uses  the new attribute {\tt FactorizedPrimes}  for number fields: it is a list containing the rational prime numbers that have already  been decomposed in $\Z_K$ already.
The running time of the  factorizations above and the local indexes are gathered in the following table:
$$
\begin{array}{|c|r|r|}
\hline
\rm{Ideal} &\rm{Index} &\rm{Time}\\
\hline\hline
2\Z_K & 26235 &0.36s \\
\hline
3\Z_K & 0& 0.61s\\
\hline
5\Z_K &20& 0.63s\\
\hline
127\Z_K & 0& 1.29s\\
\hline
313\Z_K &0& 3.69s\\
\hline
743\Z_K &0& 6.47s \\
\hline
4886229527\Z_K&0 & 6.96s \\
\hline
p\Z_K &0& 60s\\
\hline
\end{array}
$$

Thus, we need less than 90 seconds to see that the discriminant of $K$ is
$$
\mbox{Disc}(K)=
2^{1470}3^{50}5^{1960}127^{50} 313^{50} 743^{50} 4886229527^{50}p^{50}.
$$

While we cannot expect to factor the ideals $I=(w^3+50)\Z_K,$ $J=(w+10)\Z_K$ in a reasonable time, because this implies the factorization of their norm, we can factor their sum in a moment:
\magma{
> I:=ideal(K,w\pot3+50);\\
> J:=ideal(K,w+10);\\
> N:=I+J;\\
> time Factorization($\sim$N);\\
Time: 0.030\\
> N\`{}FactorizationString;\\
P(2,1)\pot2*P(2,2)\pot2*P(2,3)\pot2*P(2,4)\pot2*P(2,5)\pot2
}
This output means that all the primes {\tt P(2,1),\dots, P(2,5)} appearing in the list {\tt K`PrimeIdeals[2]} divide $I+J$ with exponent 2.  Note that, besides defining them with the new function {\tt ideal}, we manipulate the ideals with the standard {\tt Magma} notations. We can also multiply and divide fractional ideals in $K$.
\magma{
> R:=ideal(K,2);\\
> Factorization($\sim$R);\\
> time S:=R/N;\\
Time: 0.280\\
> S\`{}FactorizationString;\\
P(2,1)\pot8*P(2,2)\pot8*P(2,3)\pot36*P(2,4)\pot36*P(2,5)\pot36
}
Generators of integral ideals are efficiently computed, but the computation of generators of non-integral ideals could be quite time-costing. In the above example, it would take much time to determine generators for $N^{-1}$, but:
\magma{
> time gens:=TwoElement(S);\\
Time: 0.000
}

Our last computation in this field in this number field will be a $p$-adic square root.
%The factorization of $5\Z_K$ is:
%\magma{
%> Factorization(ideal(K,5));\newline
%[    [ 5,1, 25 ],  [ 5,2, 20 ],    [ 5,3, 5 ],    [ 5,4,\nolinebreak 20\nolinebreak ],\newline
%    [ 5,5, 5 ],    [ 5,6, 25 ],    [ 5,7, 25 ],    [ 5,8, 25 ]\nolinebreak ]
%}
We take the first two prime ideals $\p_1, \p_2$ in the list {\tt K\`{}PrimeIdeals[5]}, and check whether $w/2$ is a square in their residue fields;
\magma{
> P1:=K\`{}PrimeIdeals[5,1];\\
> P2:=K\`{}PrimeIdeals[5,2];\\
> IsSquare(w/2 mod P1);\\
false\\
> t,x0:=IsSquare(w/2 mod P2);t;x0;\\
true \\
2*z1+4\\
}
After Hensel's lemma, we now that $w/2$ has a square root in $K_{\p_2}$, but not in $K_{\p_1}$. We can compute this root with Newton's method, but we will need to lift elements from the residue field of $\p_1$ to $K$. This task is performed by the routine {\tt Lift}. For instance, to lift {\tt x0} we make:
\magma{
> X:=Lift(x0,P2);
}
\noindent
The following code computes  $\sqrt{w/2}\in K_{\p_2} (\operatorname{mod} \p_2^{10})$:
\magma{
for n in [1..10] do\\
\hphantom{dddd}  lambda:=Lift((2*X mod P2)\pot(-1),P2);\\
\hphantom{dddd}  X:=X-(X\pot2-w/2)*lambda;\\
end for;\\
print PValuation(X\pot2-w/2,P2);
}
\noindent
The whole computation takes about 60 seconds.

\subsection{A {\em fat} number field coming from the modular world}

The degree of a number field has a strong influence in the cost of the computations on it, but it is not the only handicap, nor the principal one. The discriminant turns out to be the highest obstacle in many cases, since it may be a very large integer, impossible to factor. The following example illustrates this situation.

Let $f\in \Gamma_0(1)$ be the weight 76 modular form given by the Magma command
\magma{
> f := Newforms(CuspidalSubspace(ModularForms(1,76)))[1,1];
}
It is defined over the number field $K=\Q(\theta)$, where
$\theta\in\overline{\Q}$ is a root of the polynomial
$$
\begin{array}{rl}
F(x)=x^6 +& 57080822040x^5 - 198007918566571424544768x^4 \\
&- 11405115067164354385292006554337280x^3 \\
  & + 9757628454131691442128845013041495838774263808x^2  \\
    &+290013995562379500498435975003716024800114593761580810240x\\
   &- 92217203874207784163935379997152082331434364841943058919508374716416.
\end{array}
$$
We can access this number field in Magma with the command:
\magma{
> K := Parent(Coefficient(f,1));
}

The discriminant of $F(x)$ is
$$
\operatorname{disc}(F)=
2^{264} 3^{72} 5^{16} 7^{16} 11^2 13^2 17^4 19^2 43^2 59\cdot 193^2 293\cdot 391987^2 4759427^2 137679681521^2M,
$$
where $M$ is a composite integer of 135 decimal figures which we have not been able to factorize. J. Rasmussen asked us (\cite{Rasmussen}) for a test to check certain divisibility conditions on the ring of integers of the  number field $K=\Q(\theta)$,   related to his work on congruences satisfied by the coefficients of certain modular forms.

Since we cannot factor completely the discriminant,   the {\tt Magma} function {\tt MaximalOrder}   is  useless. But we can do many computations in $K$ with our package. Let us perform a Chinese Remainder computation. We first define the number field and find the decomposition of 3 in the ring of integers:

\magma{
> Attach("dir/+ideals.m");\\
%> F:=
%x\pot6 + 57080822040*x\pot5 - 198007918566571424544768*x\pot4 \\
%\qquad- 11405115067164354385292006554337280*x\pot3 \\
%\qquad  +  9757628454131691442128845013041495838774263808*x\pot2 \\
%\qquad + 290013995562379500498435975003716024800114593761580810240*x \\
%\qquad-    92217203874207784163935379997152082331434364841943058919508374716416;\\
%> K:=NumberField(F);\\
> I:=ideal(K,3);\\
> Factorization($\sim$I);\\
> I\`{}FactorizationString;\\
P(3,1)*P(3,2)*P(3,3)*P(3,4)*P(3,5)
}

The factorization routine automatically finds generators for the prime ideals it encounters. They are stored in the attributes {\tt IntegralGenerator } and {\tt Generator} of the {\tt IdealRecord} record subjacent to the ideals, but they can be accessed with the {\tt TwoElement} function:
\magma{
> TwoElement(K\`{}PrimeIdeals[3,1]); \newline
[
    3,\\
    1/177147*(2*w\pot5 + 1815*w\pot4 + 586980*w\pot3 + 732159*w\pot2 + 658287*w + 1535274)
]
}
Consider the following Chinese Remainder problem: find  an element $\alpha\in K$ satisfying
$$
\begin{array}{lll}
\alpha\equiv w(\operatorname{mod}{\p_1}),
\qquad &\alpha\equiv w^2(\operatorname{mod}{\p_2^2}),
\quad &\alpha\equiv w^3(\operatorname{mod}{\p_3^3}),\\    \\
\alpha\equiv w^4(\operatorname{mod}{\p_4^4}),
&\alpha\equiv 1(\operatorname{mod}{\p_5}),
\end{array}
$$
where $\p_1,\dots,\p_5$ are the prime ideals above $3$ in $\Z_K$, ordered according to the output of  our factorization function. We can solve the problem using our extended version of the standard {\tt CRT} function of Magma. This extended version takes a list of numbers in a number field and a list of ideals in its ring of integers, and solves the corresponding congruences. The following piece of code finds a solution to our problem:

\magma{
>ids:=[K\`{}PrimeIdeals[3,1],K\`{}PrimeIdeals[3,2]\pot2,\\ K\`{}PrimeIdeals[3,3]\pot3,K\`{}PrimeIdeals[3,4]\pot4,K\`{}PrimeIdeals[3,5]];\\
>time r:=CRT([w,w\pot2,w\pot3,w\pot4,1],ids);r;\\
Time: 0.010    \\
1/19683*(1132037*w\pot5 + 1583958*w\pot4 + 498663*w\pot3 + 923157*w\pot2 + 649539*w + 354294)
}

We can check the result using the {\tt PValuation} routine:
\magma{
>PValuation(r-w,K\`{}PrimeIdeals[3,1]);\\
7\\
>PValuation(r-w\pot2,K\`{}PrimeIdeals[3,2]);\\
4\\
>PValuation(r-w\pot3,K\`{}PrimeIdeals[3,3]);\\
4\\
>PValuation(r-w\pot4,K\`{}PrimeIdeals[3,4]);\\
4\\
>PValuation(r-1,K\`{}PrimeIdeals[3,5]);\\
1
}

A second partial test of the result can be done reducing the element {\tt r} modulo $\p_1$ and $\p_2$:
\magma{
>(r-w) mod K`PrimeIdeals[3,1];\\
0\\
>(r-1) mod K`PrimeIdeals[3,5];\\
0
}

\subsection{A  two-parametric family of examples}\label{example3}

We now present  a two-parametric family of polynomials, which constitutes a good bench mark for our package, since the parameters involved give the chance to grade the difficulty of the computations at our convenience.

Let $\theta_n\in\overline\Q$ be a root of the polynomial $g_n=x^n+x-1$, and let $K_n=\Q(\theta_n)$. We assume that $g_n$ is irreducible. For every $r\in\Z$, the element $\theta_{n,r}:=r\theta_n$ is also a primitive element of $K_n$, with minimal polynomial $f_n^r=x^n+r^{n-1}x-r^n$. The discriminants of these polynomials are:
$$
\begin{array}{l}
\operatorname{Disc}(g_n)=(-1)^{(n-1)(n-2)/2}(n^n+(n-1)^{n-1})),\\
\operatorname{Disc}(f_n^r)=r^{n(n-1)}\operatorname{Disc}(g_n).
%\operatorname{Disc}(f_n^r)=(-1)^{(n-1)(n-2)/2}r^{n(n-1)}(n^n+(n-1)^{n-1})).
\end{array}
$$
Suppose we are given  $f_n^r$ as defining polynomial for $K_n$.
Those primes dividing $r$ but not $\Delta(g_n)$ will not ramify in $K_n$  and  we can make the local index arbitrarily large for them. Thus, an (im)proper election of $r$ can make computations in the field $K_n$ as difficult as we want.

\newpage
We will use the following function to make some tests:
\magma{
> function test(n,r)\\
function> K:=NumberField(x\pot n+r\pot(n-1)*x+r\pot n: Check:=false);\\
function> a:=Factorization(ideal(K,r));\\
function> return [K\`{}LocalIndex[p]:~p~in~PrimeDivisors(r)];\\
function> end function;
}
The function takes the degree {\tt n} and the multiplier {\tt r}, and determines the  decomposition of the primes dividing {\tt r} in the field $K_n$, defined by means of the polynomial $f_n^r$, and it returns the local index of the discriminant of $f_n^r$ at them. The following table collects the running times of the test for some values of $n$ and $r$:

$$
\begin{array}{|r||c|c|c|c|c|}
\hline
r\backslash n & 100 & 201 & 300 &  400 & 501 \\
\hline\hline
2^{50} &  0.02 & 0.05 & 0.33 & 1.04 & 0.69\\
\hline
2^{100} &  0.05 &  0.1 &  0.85 & 2.68 & 1.77\\
\hline
2^{200} &   0.07 & 0.28 & 2.23 & 7.4 & 4.73 \\
\hline
2^{400} &   0.17 & 0.76 & 6.15 & 20 & 13.71 \\
\hline
6^{50} & 0.13 & 0.77 &1.7 &2.9 & 9.77\\
\hline
6^{100} & 0.31 & 2.04 & 4.57 &7.75 & 28.09\\
\hline
6^{200} & 0.5 & 5.71 & 12.78 & 21.75 & 78.38\\
\hline
30^{50} &0.35 & 1.49& 5.02 & 6.980& 20.910\\
\hline
30^{100} &  0.82 & 3.91& 13.7& 18.34 & 58.81 \\
\hline
30^{200} & 2.24 & 10.95 & 38.49 & 49.66 & 630.81\\
\hline
210^{50} & 0.64 & 2.31 & 9.63 & 10.06  & 35.25 \\
\hline
210^{100} & 1.48 & 6.03 & 26.62 & 25.98 & 99.41 \\
\hline
2310^{100} & 2.63 & 12.25 & 31.72 & 44. 33  & 161.55\\
\hline
\end{array}
$$
For $r=m^k$ and fixed $n$, the local indices at the primes dividing $m$ are all equal to $kn(n-1)/2$. For instance, the largest local index in the table was 50100000, in the case $k=400$, $n=501$.

%always coincide; the table below shows the common local exponent of this indices:

%$$\begin{array}{|r||c|c|c|c|c|c|}\hline
%r\backslash  n & 100 & 201 & 300 &  400 & 501 \\\hline\hline
%m^{50} & 247500 & 1005000 & 2242500 & 3990000 & 6262500 \\\hline
%m^{100} & 495000 &  2010000 &  4485000  &  7980000 & 12525000\\\hline
%m^{200} &  990000 &  4020000 &  8970000 &  15960000  & 25050000\\\hline
%m^{400} &  1980000 & 8040000 & 17940000  &31920000 & 50100000 \\\hline
%\end{array}$$

We can make some comparisons with the standard Magma routines for number fields. For instance, we have compared the running time of Magma's command {\tt MaximalOrder} and our routine {\tt IntegralBasis}, with the following results\footnote{These running times do not include the time spent factoring the discriminant of the polynomial $f_n^r$; which is performed before the computation of the maximal order.}:

$$
\begin{array}{|c|c|c|c|}
\hline
n & r & {\tt MaximalOrder} & {\tt IntegralBasis} \\
\hline
5 & 2^{50} & 0.06 & 0.01 \\
\hline
5 & 2^{100} & 0.25 & 0.02\\
\hline
5 & 2^{200} & 1.77 & 0.12\\
\hline
5 & 2^{400} & 19.58 & 0.84\\
\hline
5 & 6^{50} & 0.55 & 0.04 \\
\hline
5 & 6^{100} & 2.81 & 0.23 \\
\hline
5 & 6^{200} & 23.76 & 0.01 \\
\hline
5 & 6^{400} & 270& 0.02\\
\hline
10 & 2^{50} & 1.6 & 0.18 \\
\hline
10 & 2^{100} & 14.37 & 1.11\\
\hline
10 & 2^{200} &185 & 0.1\\
\hline
20 & 2^{50} & 21.37 & 0.06\\
\hline
20 & 2^{100} & 185 & 0.1\\
\hline
40 & 210^2 & 168 & 0.71 \\
\hline
40 & 210^5 &1066 & 3.28 \\
\hline
\end{array}
$$\medskip

One may wonder whether Magma's factorization routine is faster than ours if the maximal order has already been computed. The following table compares the performance of both routines for $n=40$ and multipliers   $r_1=210^2$ and $r_2=210^5$:

$$
\begin{array}{|l|c|c|}
\hline
\multicolumn{3}{|c|} {r=210^2}\\
\hline
& {\tt Magma} & {\tt +Ideals}\\
\hline
2\Z_K & 0.35 & 0.02\\
\hline
3\Z_K & 1.77 & 0.00 \\
\hline
5\Z_K & 0.91 & 0.01 \\
\hline
7\Z_K & 1.11 & 0.01\\
\hline
\end{array}
\qquad
\begin{array}{|l|c|c|}
\hline
\multicolumn{3}{|c|} {r=210^5}\\
\hline
& {\tt Magma} & {\tt +Ideals}\\
\hline
2\Z_K & 0.36& 0.01 \\
\hline
3\Z_K & 1.78& 0.01\\
\hline
5\Z_K & 0.92& 0.00\\
\hline
7\Z_K & 1.16& 0.01\\
\hline
%25788481\Z_K & 1.26 &  0.01\\
%\hline
%p_2\Z_K & 2.95s \\
%\hline
\end{array}
$$

The local index at these primes is 2340 for $r=210^2$ and 3900 for $r=210^3$.

\subsection{The other way round: how to construct hard examples}\label{secDesign}

This last example illustrates the reverse use of the package to construct examples of number fields with prescribed decomposition of a given set of primes, determining the flux of the Montes algorithm for them. This section may be skipped for the reader not specially interested in the theory behind the package.

According to the philosophy introduced in \cite{GMN3}, a prime ideal $\p$ in a number field can be canonically represented as a list
$$
\p=[p;\phi_1,\dots,\phi_r,\phi_\p],
$$
where the $\phi_k$ are certain monic polynomials with integral coefficients, determined by Montes algorithm, when applied to the polynomial
defining the number field and the rational prime $p$. These polynomials are the main data in the  type corresponding to the ideal $\p$; the parameter $r$ is the order of this type.

Let us explain how to build a type $\ty$ with the package. We fix the prime $p=2$ and the irreducible polynomial $\phi_1=y+1\in\ff{2}[y]$, and initialize our type

\magma{
> Attach("/dir/+Ideals.m");\\
> ZX<x>:=PolynomialRing(Integers());\\
> FY<y>:= PolynomialRing(GF(2));\\
t,Y,z:=InitializeType(2,y+1);
}

The variable {\tt t} will contain our type,  the list $z$ will have the list of generators of the successive residue fields to be built, and the list $Y$ will be the list of indeterminates generating the corresponding polynomial rings.
At this moment {\tt t} is the order one type given by the lift $\phi_1=x+1\in\Z[x]$ of $y+1$. We proceed to enlarge it, constructing a polynomial $\phi_2\in\Z[x]$ whose $\phi_1$-Newton polygon has a unique side with slope 1/2 and residual polynomial $y^2+y+1$:

\magma{
> EnlargeType(1,2,y\pot2+y+1,~t,~Y,~z);
}
The polynomial $\phi_2$ is
\magma{
> t[2]\`{}Phi;\\
x\pot4 + 4*x\pot3 + 8*x\pot2 + 8*x + 7
}

We go one step further, and enlarge the type with a polynomial $\phi_3\in\Z[x]$ having a one-sided $\phi_2$-Newton polygon, with slope 2/3 and residual polynomial  $Y[2]^2+z[2]Y[2]+z[1]$, where $\ff{\ty}=\ff{2}(z[2])=\ff{2}[y]/(y^2+y+1)$, and $z[1]=1$.

\magma{
> EnlargeType(2,3,Y[2]\pot2+z[2]*Y[2]+z[1],~t,~Y,~z);\\
> t[3]\`{}Phi;\\
x\pot24 + 24*x\pot23 + 288*x\pot22 + 2288*x\pot21 + 13482*x\pot20 + 62664*x\pot19 + 238736*x\pot18 + 765072*x\pot17 + 2100447*x\pot16 + 5005808*x\pot15 + 10455232*x\pot14 +
    19267808*x\pot13 + 31474060*x\pot12 + 45694864*x\pot11 + 59023136*x\pot10 + 67789216*x\pot9 + 69069935*x\pot8 + 62166904*x\pot7 + 49095584*x\pot6 + 33671984*x\pot5 +
    19743658*x\pot4 + 9650728*x\pot3 + 3776208*x\pot2 + 1092816*x + 208081
}

Let us build a second type:
\magma{
> tt,YY,zz:=InitializeType(2,y\pot2+y+1);       \\
> EnlargeType(3,2,y+1,$\sim$tt,$\sim$YY,$\sim$zz);\\
> EnlargeType(3,1,YY[2]\pot2+YY[2]+zz[1],$\sim$tt,$\sim$YY,$\sim$zz);\\
> EnlargeType(1,3,YY[2]+zz[2],$\sim$tt,$\sim$YY,$\sim$zz);\\
> tt[4]\`{}Phi;     \\
x\pot24 + 12*x\pot23 + 102*x\pot22 + 616*x\pot21 + 3045*x\pot20 + \dots
%12432*x\pot19 + 43866*x\pot18 +
%    136068*x\pot17 + 379554*x\pot16 + 975076*x\pot15 + 2323662*x\pot14 + 5215344*x\pot13 +
%    10861517*x\pot12 + 21141216*x\pot11 + 38275254*x\pot10 + 64921532*x\pot9 + 105639186*x\pot8
%    + 156993756*x\pot7 + 215634914*x\pot6 + 248842656*x\pot5 + 286107093*x\pot4 +
%    430926088*x\pot3 + 264711118*x\pot2 + 111773236*x + 3581577
}
We can combine any number of types to produce a polynomial whose family of attached types will be the given ones.
\magma{
> pol:=CombineTypes([t,tt]);\\
> pol;\\
x\pot48 + 36*x\pot47 + 678*x\pot46 + 8808*x\pot45 +\dots
% 88143*x\pot44 + 720744*x\pot43 + 4994470*x\pot42
%    + 30062772*x\pot41 + 160009107*x\pot40 + 763367360*x\pot39 + 3299806876*x\pot38 +
%    13040321024*x\pot37 + 47472259240*x\pot36 + 160260265616*x\pot35 + 504644562544*x\pot34
%    + 1489881733576*x\pot33 + 4142471925501*x\pot32 + 10887958805676*x\pot31 +
%    27136789803954*x\pot30 + 64294052124032*x\pot29 + 145084943565289*x\pot28 +
%    312287280202064*x\pot27 + 641861370000774*x\pot26 + 1260691506341492*x\pot25 +
%    2367236108894694*x\pot24 + 4250016133967972*x\pot23 + 7294378384412750*x\pot22 +
%    11963754692503864*x\pot21 + 18740061736377609*x\pot20 + 28011428060782528*x\pot19 +
%    39906029090031338*x\pot18 + 54090423271504564*x\pot17 + 69581543299059973*x\pot16 +
%    84660397604187616*x\pot15 + 97002579970349376*x\pot14 + 104114133096457160*x\pot13 +
%    104042743523340888*x\pot12 + 96144295886586216*x\pot11 + 81537310364202300*x\pot10 +
%    62925486917088024*x\pot9 + 43760401059255179*x\pot8 + 27100773748847660*x\pot7 +
%    14721791132600574*x\pot6 + 6873774725531288*x\pot5 + 2679467144987359*x\pot4 +
%    835591889852560*x\pot3 + 190753714537150*x\pot2 + 27171891370948*x +
%    2582249878086908589655919172003011874329705792829223512830659356540647622016841194629645353280137831435903172718005617113
}
You can check that the polynomial {\tt pol} has the two types $\ty, \ty\ty$ at the prime $p=2$, and follow the flux of the algorithm using the flag {\tt montestalk}  to set the level of verbosity of the package:
\magma{
 >    SetVerbose("montestalk",4);\\
 >    K<w>:=NumberField(pol);\\
 >    Factorization(ideal(K,2));\\
Analyzing irreducible factor modulo p:  Y0 + 1\\
++++++++++++++++++++++++++++++++++++++++++++++++\\
++++++++++++++++++++++++++++++++++++++++++++++++\\
Analyzing type of order  1\\
Phir= x + 1\\
Sides of Newton polygon: [\newline
    [ -1/2, 0, 12, 24, 0 ]\newline
]\newline
----------------------\\
Analyzing side  [ -1/2, 0, 12, 24, 0 ]\\
Slope:  -1/2\\
Origin: ( 0 , 12 )\\
End point: ( 24 , 0 )\\
----------------------
Monic Residual Polynomial= Y0\pot12 + Y0\pot10 + Y0\pot6 + Y0\pot2 + 1     \\
Factors of R.P.:= [\newline
    <Y0\pot2 + Y0 + 1, 6>\newline
]\\
\dots
}

In some sense, these constructions lead to number fields specially ill posed for Montes algorithm, since the order of the types dominates the complexity of the computations. The relation:
$$
\deg(K)=\sum_{\ty\in \mathbf{T}(f,p)}\prod_{k=0}^{\sharp\ty} e_k^{\ty}f_k^{\ty},
$$
between the degree of a number field and certain entries $e$ and $f$ (partial ramification indices and residual degrees) of the types attached to the prime ideals lying over the prime number $p$, shows that the growth of the  order $\sharp\ty$ of these types is logarithmic with respect to the degree.

\section{Description of the package}

\subsection{New data types and attributes for ideals and number fields}

The centre of the attention of the package are number fields and its fractional ideals. Since computations tend to be quite time-consuming, the important information obtained from any operation is saved. This implies the introduction of new attributes for number fields, and new data types for handling their fractional ideals.

The new attributes for number fields are:
\begin{itemize}
\item \texttt{FactorizedPrimes}: List containing the rational primes already decomposed in~$K$.
\item \texttt{PrimeIdeals}: An associative array, with one entry per each rational prime $p$ in the {\tt FactorizedPrimes}; it  contains the prime ideals in $\Z_K$ dividing $p$.
\item \texttt{LocalIndex}: Associative array which contains, for every rational prime in {\tt FactorizedPrimes}, the local exponent of the index of the polynomial defining the number field.
\item \texttt{Discriminant} The discriminant of the number field.
\item \texttt{FactorizedDiscriminant}: The factorization of the discriminant of the number field.
\item \texttt{pIntegralBasis}: An associative array, which for every rational prime $p$ in {\tt FactorizedPrimes} stores a $p$-integral basis of the maximal order of the field.
\item \texttt{IntegralBasis}: A list with a basis of the maximal order of the fields.
\end{itemize}
These attributes are assigned and expanded when the corresponding objects are computed. For instance, an integral basis of $K$ is not computed unless the function {\tt IntegralBasis} is used, but then the last four attributes are set, and they are never recomputed.

The fractional ideals are handled by means of two different record structures:
the {\tt PrimeIdealRecord} and the {\tt IdealRecord}.  The first one is designed to contain the essential information produced by Montes algorithm for every prime ideal, while the second one is used for general ideals. Their concrete
fields are:

\magmad{
PrimeIdealRecord:=recformat<\\
Parent:  //   Number field where the prime ideal lies.\\
Pol:   // Polynomial defining the number field.     \\
Prime: // Rational prime $p$ below the prime ideal.\\
e: //  Ramification index of the prime ideal.\\
f: // Residual degree of the prime ideal.\\
Generator: // The ideal is generated by Prime and this element.\\
%bp: //  Technical entry, for internal use only.\\
Position:  Position of the prime ideal in the list\\ \hphantom{dddddddd} K`PrimeIdeals[$p$].\\
Type:  // Type corresponding to the prime ideal.\\
>;
}
\magmad{
IdealRecord:=recformat<\\
    Parent: Number field where  the  ideal lies.\\
    Generators:  // Generators given in the definition of the ideal.\\
    Norm:  // Norm of the ideal.\\
    Radical: // List of rational primes dividing the norm \\ \hphantom{dddddddd} of the ideal.\\
    IsIntegral: // Flag to describe the integrality of the ideal.\\
    IsPrime: // Flag to describe the primality of the ideal.\\
    Factorization: // List with the prime ideals dividing the \\ \hphantom{dddddddd}  ideal, and the corresponding exponents.\\
    FactorizationString:  // Factorization of the ideal in \\ \hphantom{dddddddd} pretty format.\\
    IntegerGenerator:  // The ideal is generated by it and the \\ \hphantom{dddddddd} element in the next field.\\
    Generator:  // The ideal is generated by it and the element \\ \hphantom{dddddddd} in the previous field.\\
>;
}

The fields in these structures are assigned only when its computation is trivial or  strictly necessary.
Once  ideals or  prime ideals are created, they can be manipulated using the standard notation of Magma for ideal arithmetic.

The trivial ideals $(0)$ and $\Z_K$: an ideal is the zero ideal if its {\tt Generators} field is assigned, and all its elements are zero; and an ideal is the full ring of integers if its {\tt Factorization} field is assigned and is an empty list. One must check whether these fields are assigned, since the ideal may be the result of an operation between other ideals, which does not assign any of these fields. The functions {\tt IsZero} tests whether an ideal is zero.

\subsection{Main functions}

The  essential tasks of the package is the creation and factorization of fractional ideals in number fields;  they are performed by means of the functions {\tt ideal} and {\tt Factorization}.

The purpose of the function {\tt ideal} is the creation of a general fractional ideal, given a generator or a list of generators. Accordingly, it admits the two forms:
\magma{
ideal(K, alpha1 );\\
ideal(K, [alpha1,\dots,alphar]);
}
\noindent
where {\tt K} is a number field and {\tt alpha1,\dots,alphar} are elements of {\tt K}. The function returns an {\tt IdealRecord} structure, with the fields {\tt Parent} and {\tt Generators} assigned.

Once we have created a fractional ideal {\tt I}, we can call the extended version of {\tt Factorization} to find its
decomposition into prime ideals. We can pass the ideal as a direct argument to the function or as a reference argument:
\magma{
Factorization(I);\\
Factorization($\sim$I);}
They do the same computations, but the procedural form assigns the fields {\tt Factorization}, {\tt FactorizationString} and {\tt IsPrime} in the record given by {\tt I} (if they are not yet assigned) and returns nothing, while the first form returns the value of the field {\tt Factorization}. The behaviour of the functions depend on the data known for the ideal:
If the field {\tt Factorization} is assigned, the function returns its value.
Otherwise, the greatest common divisors of the norms of the generators of the ideal is computed and factored. Every prime number dividing this gcd is decomposed in the number field (if it is not already decomposed) and the valuation of the given ideal at the prime ideals above these prime numbers are computed with the {\tt PValuation} function.

The decomposition of a rational prime in a number field is performed by the {\tt Montes} function; the command
\magma{
 Montes(p,K);}
\noindent applies Montes algorithm to the polynomial defining the number field {\tt K} and the rational prime {\tt p}, to find the prime ideal decomposition of {\tt p} in {\tt K}. It appends {\tt p} to the list {\tt K\`{}FactorizedPrimes}, and fills the entries corresponding to {\tt p} in  the associative arrays {\tt K\`{}LocalIndex}, and {\tt K\`{}PrimeIdeals}.

A key tool for the factorization is the {\tt PValuation} function. It accepts an element of a number field and a prime ideal in the field (given as a {\tt PrimeIdealRecord} or an {\tt IdealRecord}) and returns the valuation of the ideal generated by the element at the given prime.

Since the discriminant of a number field is a very important invariant, the package includes the function {\tt TrueDiscriminant} to compute it. Note that it has to factor the  discriminant of the polynomial defining the number field, which may be quite time-consuming or even impossible.

\subsection{Ideal arithmetic}
We can add (+), multiply (*) or divide (/) two ideals and raise a non-zero ideal to any integral power (\pot) using the standard operators.  The mantra for the design of the arithmetic of ideals is {\em "don't waste time in unnecessary hard computations"}. Thus, for instance, the addition of two ideals simply joins its list of generators and only if the factorization of both ideals is known also finds the factorization of the result. The multiplication of two ideals combines the factorizations of the ideals (computing them if necessary), but does not find generators for the product.

The ideals can also be compared in the usual way: the operators {\tt eq} and {\tt subset} have extended versions for such task. They simply compare the factorizations of the ideals to be compared. The norm of an ideal and  the rational primes dividing it are given by the functions {\tt Norm}, and {\tt RationalRadical}.

\subsection{Generators}

It is well known that any fractional ideal in a number field can be generated by a rational integer and a second integral element in the field. The usual Magma  function  {\tt TwoElement} has been extended to apply to our new ideals. It admits two forms
\magma{
TwoElement(I);\\
TwoElement($\sim$I);
}
They perform the same computation, but the procedural version also assigns the fields {\tt Generators}, {\tt IntegralGenerator} and {\tt Generator} of the given ideal. The factorization of the ideal is computed if it is not known.

\subsection{Congruences. Chinese remainder theorem}

The package provides facilities for congruential computations for prime ideals. Essentially, the package provides the reduction map $\Z_K\rightarrow \Z_K/\p$ from the maximal order of the number field $K$ to the residue field of any prime ideal $\p$. The residue field of a prime ideal is accessed through the command:
\magma{
ResidueField(P);}

The reduction mod $\p$ of a $\p$-integral element  $\alpha\in K$ is computed with
\magma{
alpha mod P;}

An element $z\in\ff{\p}$ in the residue field of $\p$ may be lifted to an integral element in $\Z_K$ with the command
\magma{
Lift(z,P);}

The Magma {\tt ChineseRemainderTheorem} function is extended in the package to work with the new ideal records. Its syntax is
\magma{
CRT([alpha1,\dots,alphar],[id1,\dots,idr]);
}
\noindent
where {\tt alpha1,\dots,alphar} are integral elements of a number field and {\tt id1,\dots, idr} are integral ideals (given as {\tt IdealRecord}'s or as {\tt PrimeIdealRecord}'s). The function returns an element {\tt alpha} in the same number field, such that {\tt alpha-alphak} is divisible by {\tt idk}. Note that the function does not test the integrality of the given elements.

\subsection{Integral Basis}

A conjectural algorithm for the computation of integral bases of number fields can be derived from Montes algorithm. It has been developed in \cite{GMNbasis}. The algorithm always returns a correct answer, since it checks its validity; in case of fail, the user is warned. Extensive tests of the algorithm have been done, without no counterexample found.

The  functions included in the package for this task are:
\magma{
pIntegralBasis(K,p);\\
SIntegralBasis(K,[p1,\dots,pr]);\\
IntegralBasis(K);\\
}
The basic one is {\tt pIntegralBasis}, which takes a number field {\tt K} and a prime number {\tt p} and computes a {\tt p}-integral basis of the maximal order $\Z_{\tt K}$; the basis is saved in the corresponding entry of the associative array {\tt K\`{}pIntegralBasis}. The prime {\tt p} is decomposed in {\tt K}  only if this has not been done yet.

 The {\tt SIntegralBasis} functions takes a list of primes as second input, and returns a local integral basis for this set of primes. It calls {\tt pIntegralBasis} for every given prime, and glues the local bases with a Chinese Remainder computation.

Finally, the function {\tt IntegralBasis} computes a $\Z$-basis of the maximal order of the given number field. It simply calls {\tt SIntegralBasis} with the list of primes dividing the discriminant of the polynomial defining the number field.

\subsection{Construction of {\em hard} number fields}

In order to check the main functions of the package we have developed a number of functions to create examples of number fields with prescribed behaviour  of the Montes algorithm. They allow to create and combines {\it types} from its essential data. The basic functions for this task are:
\magma{
InitializeType(p,psi);\\
EnlargeType(  h,e,psi,~t, ~Y,~z);
}
 A detailed example of its use has been presented in section \ref{secDesign}. Remind that we can use the function {\tt CombineTypes} to combine different types with the same rational prime. The function

  \magma{
  CombinePolynomialsWithDifferentPrimes(f1,p1,f2,p2,k)
 }
\noindent
 takes two polynomials {\tt f1, f2} of the same degree and two prime numbers {\tt p1, p2}, and finds a polynomial $f$ of the same common degree such that ${\tt f}\equiv {\tt f}(\operatorname{mod}\, {\tt pi}^k)$. It is useful for the construction of examples of number fields with prescribed types at different primes.

 Random examples can be easily produced using the functions:
 \magma{
 CreateRandomType(p,r);  \\
 CreateRandomMultipleTypePolynomial(p,k,r,s);
 }
 \noindent
 The first one generates a random type of order {\tt r} for the prime number {\tt p},  while the second one computes a random polynomial which has {\tt k} different types of order {\tt r} for the prime {\tt p}. This polynomial is obtained as the product of all the last $\phi$-polynomials of the randomly created types, plus an integral multiple of {\tt p\pot s} to ensure its irreducibility.

 Since the examples produced by these functions use to be enormous, some functions to write them in a format proper for Magma or TeX are included. They are {\tt ExpandTeX, ExpandPhiTeX, ExpandAllPhiTeX} and  {\tt ExpandMagma}.
Their self-explanatory name makes unnecessary any further comment.

\end{document}